\documentclass[11pt]{article}
\usepackage{amssymb,amsfonts,amsmath,latexsym,graphicx}
\newtheorem{theorem}{Theorem}[section]
\newtheorem{proposition}[theorem]{Proposition}

\newtheorem{corollary}[theorem]{Corollary}

\newcommand{\be}{\begin{equation}}
\newcommand{\ee}{\end{equation}}
\newcommand{\bea}{\begin{eqnarray}}
\newcommand{\eea}{\end{eqnarray}}
\newcommand{\la}{\label}
\def\R{\mathbb{R}}
\def\N{\mathbb{N}}

\def\a{\alpha}
\def\xa{\alpha}
\def\xs{\sigma}
\def\xd{\delta}
\def\xl{\lambda}
\def\xe{\epsilon}
\def\xg{\gamma}
\def\xO{\Omega}
\def\xb{\beta}
\def\l{\lambda}
\def\pa{\partial}
\def\finedim{{\unskip\nobreak\hfil\penalty50
\hskip2em\hbox{}\nobreak\hfil\hbox{{ \vrule height.9ex   width.8ex
depth-.1ex  }\qquad} \parfillskip=0pt
\finalhyphendemerits=0\par\medskip}}

\title{On a class of  weighted anisotropic Sobolev
inequalities}

\author{\Large Stathis Filippas$^{1,4}$ \& Luisa Moschini$^{2}$
\& Achilles Tertikas$^{3,4}$  \\
                                                                           \\
        Department of Applied Mathematics$^{1}$ \\
         University of Crete,
         71409 Heraklion,  Greece \\
        filippas@tem.uoc.gr\\
                                          \\
                                          Dipartimento di Metodi e Modelli Matematici''$^2$,\\
 University of Rome ``La Sapienza'', 00185 Rome, Italy \\
moschini@dmmm.uniroma1.it \\
\\
 Department of Mathematics$^{3}$ \\
         University of Crete,
         71409 Heraklion,  Greece \\
          tertikas@math.uoc.gr\\
                                     \\
        Institute of Applied and Computational Mathematics$^4$, \\
        FORTH, 71110 Heraklion, Greece \\
    \\ }

\begin{document}
\date{\today}
\maketitle

\begin{section}{Introduction and main results}

In this article, motivated by  the  work of Caffarelli and Cordoba  [CC] in phase
transitions analysis, we prove new  weighted {\it
anisotropic } Sobolev type inequalities, that  is Sobolev type
inequalities where  different derivatives have different weight
functions.

Phase transitions or interfaces appear in physical problems when
two different states coexist and there is a balance between two
opposite tendencies: a diffusive effect that tends to mix the
materials and a mechanism that drives them into their pure state,
which is typically given by a nonnegative potential $F(x,u)$,
denoting the energy density of the configuration $u$. For example
it is known that minimizers of the functionals
$$J_\epsilon(u):=\int_{\Omega} \{\epsilon^2 |\nabla u|^2 +F(x,u)\}
dx,$$ for $0<\epsilon<1$, $F(x,u)=(1-u^2)^\delta_{+}$, and
$\Omega\subset \R^N$ open and bounded, develop free boundaries
 if $0<\delta < 2$, while generate exponential convergence
to the states $\pm 1$ if $\delta=2$, that is in the case connected
to the Ginzburg-Landau equation, see   \cite{CC}.

The main  results of  \cite{CC}  are concerned with the study of
regularity properties  of interfaces. Their results are closely
related to a conjecture of De Giorgi according to which bounded
solutions of the Ginzburg-Landau scalar equation on the whole
space $\R^N$ that are monotone in one direction, are  one
dimensional (see \cite{DG}); in particular they concern the
question of De Giorgi under the additional assumption that the
level sets are the graphs of an equi-Lipschitz family of functions
(see \cite{MM} for the case $N=2$, see also \cite{BBG} for the
general case). In establishing these results a  central role is
played by  various anisotropic  Sobolev type inequalities, see
Propositions 4--5  in  \cite{CC}.

Moreover, the weighted anisotropic Sobolev inequalities we are
dealing with, are also intimately  connected to Sobolev
inequalities for Grushin type operators. Unweighted local version
of this type of inequalities
 have  been studied in \cite{FL1},
\cite{FL2}, as well as in \cite{FGW} where Muckenhoupt weights
were considered.

As a further motivation to the present study, we mention
 that Sobolev inequalities,    are used in the  proof of Liouville type
theorems for the corresponding linear
elliptic operators in divergence form.

For other type of  anisotropic  Sobolev type  inequalities
we refer to \cite{Ba}, \cite{Be}, \cite{Mo}.

To state our  results  let us   first introduce some
notation. We define the infinite cylinder $\mathcal H_1$ as well as the
finite cylinder $\mathcal C_{1}$  by
\bea
\mathcal H_1  & : =  &
\{(x',\lambda)\in \R^{N-1}\times \R: |x'|<1\} ,
 \nonumber \\
\mathcal C_{1}  &  :=  & \{(x',\lambda)\in \R^{N-1}\times \R:
|x'|<1,~ |\lambda|<1\}.     \nonumber
\eea
We  will  prove weighted
Sobolev inequalities on the finite cylinder $\mathcal C_1$, the
weight being a positive power of the distance function to the top
or the bottom of the cylinder $\{\l=\pm 1\}$.

Our first result is the following

\begin{theorem}\label{th1}
 Let $N\ge 2$,
 $\alpha >-1$ and  $\sigma\in (-2\alpha,2)$.   Then, for any  $Q$ with
\be\la{qcr}
 2 \leq  Q \leq Q_{cr}(N,\alpha,\sigma):=
\frac{2 \left( N + \frac{2 \a +\xs}{2-\xs}\right)}{N + \frac{2 \a
+\xs}{2-\xs} -2}, \ee
 there exists a positive constant $C=C(Q,N,\alpha,\sigma)$, such that
 for any function $f\in C^\infty_0(\mathcal H_1)$ there holds
\begin{equation}\la{11}
 \left(\int_{\mathcal C_1} (1-|\lambda|)^\alpha
|f(x',\lambda)|^Q dx' d\lambda \right)^{\frac{2}{Q}} \le C
\int_{\mathcal C_1}  (1-|\lambda|)^\alpha \left(|\nabla_{x'} f|^2
+(1-|\lambda|)^{\sigma}|\partial_\lambda f|^2\right) dx' d\lambda.
 \end{equation}
 In the limit case where $\sigma =2$,  estimate
(\ref{11})  holds for $Q=2$ and any $f\in
C^\infty_0(\mathcal H_1)$ but  fails for $Q>2$
 and $f\in C^\infty_0(\mathcal C_1)$.
\end{theorem}

In the case  $\sigma  \geq 2$  we can still have similar inequalities
when  $\alpha<-1$. More precisely   when $\xs=2$ we have

\begin{theorem} \label{th4444} Let $N\ge 2$ and $\alpha<-1$. For any $Q$ with
 $2\le Q \le \frac{2N}{N-2}$, in case $N \geq 3$, or $Q \geq 2$ in case $N=2$,
there exists a positive constant  $C=C(N,\alpha,Q)$, such that for
any function $f\in C^\infty_0(\mathcal C_1)$ there holds
\begin{equation}\label{vivi}
\left(\int_{\mathcal C_1} (1-|\lambda|)^{\alpha}
 |f(x',\lambda)|^Q dx' d\lambda
\right)^{\frac{2}{Q}} \le C   \int_{\mathcal C_1}
(1-|\lambda|)^\alpha   \left(|\nabla_{x'} f|^2
+(1-|\lambda|)^{2}|\partial_\lambda f|^2\right) dx' d\lambda .
\end{equation}
\end{theorem}

When  $\xs>2$  we obtain the same inequality but
this time
for  exponents  $Q$  that satisfy $Q \geq Q_{cr}$
 as defined in  (\ref{qcr}).  Thus, we have

\begin{theorem}\label{th1.3}
 Let $N\ge 2$,
 $\alpha <-1$ and  $\sigma\in (2, -2\alpha)$.  Then,
for any $Q$ with $Q_{cr} \leq  Q$ if $N=2$ or
 $ Q_{cr} \leq Q \leq  \frac{2N}{N-2}$ if $N \geq 3$,
 there exists a positive constant $C=C(N,Q,\alpha,\sigma)$, such that
 for any function $f\in C^\infty_0(\mathcal C_1)$ there holds
\begin{equation}\la{11.3}
 \left(\int_{\mathcal C_1} (1-|\lambda|)^\alpha
|f(x',\lambda)|^{Q} dx' d\lambda \right)^{\frac{2}{Q}} \le C
\int_{\mathcal C_1}  (1-|\lambda|)^\alpha \left(|\nabla_{x'} f|^2
+(1-|\lambda|)^{\sigma}|\partial_\lambda f|^2\right) dx' d\lambda.
 \end{equation}
\end{theorem}

When $\xa>-1$ then $(1 -|\xl|)^{\xa}$ is an $L^1(-1,1)$ function  and using
Holder's inequality one can obtain  the inequality for any $Q$ with $2 \leq  Q \leq Q_{cr}$
once it is true for $Q_{cr}$. However this is not the case when $\xa <-1$.

We  note that for $Q=2$  inequality (\ref{11.3}) is still valid as
one can see using Poincar\'e inequality in the $x'$--variables.
The validity or not of  (\ref{11.3})  for $2 < Q < Q_{cr}$ remains
an open question.

Finally,  as $\xs >2$ approaches 2,  $Q_{cr}$  approaches 2
and therefore the $Q$--interval  of validity of (\ref{11.3})
approaches the interval $[2, \frac{2N}{N-2}]$
in
complete agreement with the result of Theorem \ref{th4444}.

A central role in  the proof of the previous results,
is played by various
 weighted isotropic Sobolev inequalities in
the upper half space $\R^N_{+}:=\{(x',x_N)\in \R^{N-1}\times \R:
x_N>0\}$,   which are of independent interest. We present such a result:
\begin{theorem}\label{th2} Let either
\begin{equation}\label{cond2}
N=2,  ~~~~~~~ 2  \le Q,   ~~~~~~~~{\rm and}~~~~~~ B=A- \frac{2}{Q}
\ ,
\end{equation}
or else,
\begin{equation}\label{cond1}
N\geq 3, ~~~~~~~ 2  \le Q  \le \frac{2N}{N-2}, ~~~~~~{\rm
and}~~~~~~ B=A-1+\frac{Q-2}{2Q}N .
\end{equation}
If  $BQ+2A \neq 0$, or if $A=B=0$ then  \\
(i)There exists a positive constant $C=C(A,Q,N)$, such that for
any function $f\in C^\infty_0(\R^N_{+})$ there holds
\begin{equation}\label{due}
 \left(\int_{\R^N_{+}} x_N^{BQ}
|f(x',x_N)|^Q dx' dx_N \right)^{\frac{2}{Q}} \le C \int_{\R^N_{+}}
x_N^{2A} \left(|\nabla_{x'} f|^2 +|\partial_{x_N} f|^2  \right) \
dx' dx_N \ .
\end{equation}
(ii)  If moreover  $BQ+2A>0$, or if $A=B=0$ inequality
(\ref{due}) still holds even if  $f\in
C^\infty_0(\R^N)$.
\end{theorem}

The exponent $Q=Q(A,B,N)$ given by conditions (\ref{cond1}) and
(\ref{cond2}) is the best possible, as one can easily see arguing
by scaling $x'=Ry'$, $x_N=Ry_N$. In case $ N \geq 3$, part  (i) of
the Theorem  \ref{th2} is due to Maz'ya, see  \cite{M}, section
2.1.6. Here we will provide a simpler proof along the lines of
\cite{FMaT1}, \cite{FMaT2}, \cite{FMT}. A particular case of
(\ref{due})  has been obtained in \cite{C} under an additional
assumption on $f$,
 by different methods.

We next present a  direct consequence of Theorem \ref{th1}.

\begin{corollary} \label{th3}
For  $N\ge 2$,  $m>-1$ and $\epsilon\in(0,\frac12)$ we set
\[
C_{1,\epsilon}:=\{(x',\lambda)\in \R^{N-1}\times \R: |x'|<1,
|\lambda|<1-\epsilon^{1+m}\}\ .
\]
Let $\alpha>-1$ and  $\beta >0$ satisfy
$$-2\alpha (1+m)<\beta m<2(1+m),$$
and
\[
2 \leq P \leq P_{cr}(N,m,\alpha, \beta)
 :=  \frac{2 \left( N + \frac{2 \alpha (1+m)+\beta m}{2 (1+m)-\beta m}\right)  }
{N + \frac{2 \alpha (1+m)+\beta m}{2 (1+m)-\beta m} -2}.
\]
Then, there exists a positive constant $C=C(N,P,m,\alpha, \beta)$
independent of $\epsilon$, such that
 for any function $f\in C^\infty_0(C_{1,\epsilon})$ there holds
$$\left(\int_{C_{1,\epsilon}}
(1-|\lambda|)^\alpha |f(x',\lambda)|^P dx' d\lambda
\right)^{\frac{2}{P}} \le C
\int_{C_{1,\epsilon}}(1-|\lambda|)^\alpha \left(|\nabla_{x'} f|^2
+\frac{(1-|\lambda|)^\beta}{\epsilon^\beta} |\partial_\lambda
f|^2\right) dx' d\lambda. $$
\end{corollary}

The above corollary is in the same spirit as the results in
\cite{CC}. Indeed, when $\alpha=1$ and $\beta=2$, Corollary
\ref{th3} entails the weighted Sobolev inequality of Proposition 5
of \cite{CC} providing a precise range for the Sobolev exponent.
Analogous results can be easily obtained in case $\xa <-1$, by
using Theorems \ref{th4444} and \ref{th1.3}.

 We   next  consider  the more general case of  weighted anisotropic
inequalities  where the distance is taken from a higher  codimension
boundary. More precisely,
 for $x \in \R^N$ we write $x=(x', \l)$,
with $x' \in \R^{N-k}$ and $\l \in \R^k$,  with $1<k<N$.
Let $\xO \subset \R^k$ be a smooth bounded domain and
 $B_1=\{x': |x'|<1 \}$
be the unit ball in $\R^{N-k}$. We also set  $d=d(\l)={\rm
dist}(\l, \pa \xO)$.  In this case our main result reads

 \begin{theorem}\la{th9}
 Let $N \geq 3$, $1 < k <N$,  $\xa>-1$ and  $\xs \in (- 2 \xa,
 2)$ with $2 \xa + \xs k \geq 0$.
  Then, for any $Q$,
$$2 \leq Q \leq Q_{cr}^k
:=\frac{2(N +\frac{ 2 \xa + \xs k}{2-\xs})}{N +\frac{ 2 \xa +
\xs k}{2-\xs}-2}, $$  there  exists a positive constant
$C=C(Q,N,\xa,\xs,k)$, such that for any function  $f \in
C_0^{\infty}(B_1 \times \Omega)$ there holds
\begin{equation}\la{4.30}
\left(\int_{B_1 \times \xO} d^\alpha |f|^Q dx
\right)^{\frac{2}{Q}} \le C \int_{B_1 \times \xO} d^{\xa}
\left(|\nabla_{x'} f|^2 +d^{\sigma}|\nabla_{\l} f|^2\right) dx.
\end{equation}
\end{theorem}
The  limit  case $k=N$,   corresponds to
the following  isotropic weighted inequality
\[
\left(\int_{\xO} d^\alpha|f|^Q d\lambda\right)^{\frac{2}{Q}} \le C
\int_{\xO} d^{\xa+\sigma}|\nabla_{\l} f|^2 d\lambda,
\]
which is true  when $\xa+\xs <1$
but not  when $\xa +\xs \geq 1$; see Remark after the proof of
Theorem \ref{th9} for details.

To prove the above Theorem an important role is played by
the following
 weighted  anisotropic  Sobolev inequality in the
upper half space $\R^N_{+}$.
To state the result we  first introduce some notation.
 For $x \in
\R^N_+$, $1< k < N$, we write $x=  (x', \xl)= (x',x_N, y)$, with $x' \in
\R^{N-k}$, $x_N \in [0, \infty)$, and $y \in \R^{k-1}$. We also
write  $dx$ for $dx' d \xl= dx' dx_N dy$.

\begin{theorem}\label{th40} Let  $\xg \in \R$, and either
\begin{equation} \label{maybb}
N=2, ~~~~~~~ Q  \ge 2, ~~~~~~{\rm and}~~~~~~
B=A-1+\frac{Q-2}{2Q}(2+\xg(k-1)),
\end{equation}
or else,
\begin{equation} \label{may}
N\geq 3, ~~~~~~~ 2  \le Q  \le \frac{2N}{N-2}, ~~~~~~{\rm
and}~~~~~~ B=A-1+\frac{Q-2}{2Q}(N+\xg(k-1)).
\end{equation}
If  $BQ+2A \neq 0$ then  \\
(i)There exists a positive constant $C=C(A,Q,N, k, \gamma)$, such
that for any function $f\in C^\infty_0(\R^N_{+})$ there holds
\begin{equation}\label{due4}
 \left(\int_{\R^N_{+}} x_N^{BQ}
|f(x)|^Q dx \right)^{\frac{2}{Q}} \le C \int_{\R^N_{+}} x_N^{2A}
\left(|\nabla_{x',x_N} f|^2 + x_N^{2 \xg}
 |\nabla_{y} f|^2 \right) \ dx\ .
\end{equation}
(ii)  If moreover  $BQ+2A>0$, inequality (\ref{due4}) still holds
even if  $f\in C^\infty_0(\R^N)$.\\
\end{theorem}
We note  that the exponent $Q=Q(A,B,N,\gamma,k)$ given by
(\ref{may}) is the best possible  as one can easily check
 using the natural scaling    $x'=Rz'$,
$x_N=Rz_N$ and $y=R^{\gamma+1}w$.

Inequality (\ref{due4}) is a weighted Sobolev inequality for
Grushin type operators $\mathcal
L_\gamma:=\Delta_{x',x_N}+x_N^{2\gamma}\Delta_y$ having associated
gradient $\nabla_\gamma:=(\nabla_{x'},
\partial_{x_N}, x_N^\gamma \nabla_y)$, so that
$$|\nabla_{\gamma}g|^2=|\nabla_{x',x_N}g|^2 +x_N^{2\gamma}
|\nabla_y g|^2 \ .$$ When $\gamma\in \N$  then
$\mathcal L_{\gamma}:=\frac{\partial^2}{\partial x_1^2}
+x_1^{2\gamma} \frac{\partial^2}{\partial x_2^2}$ belongs to the
class of differential operators considered by \cite{B}; in
particular, it is hypoelliptic and satisfies a Harnack inequality
since the Lie algebra generated by the vector fields
$\frac{\partial}{\partial x_1}$ and
$x_1^{\gamma}\frac{\partial}{\partial x_2}$ has rank two at  any
point of the plane. On the other hand  when $\gamma>0$ and $-1<2A<1$,
the weight is a Muckehoupt weight and the  local version of
inequality (\ref{due4}) was considered in \cite{FGW}.
Our method  has the advantage of allowing a bigger range of values for
the parameter $A$,  in particular allowing weights outside the
Muckehoupt classes.

We finally note that  weighted Sobolev type inequalities
of the kind we present in this work, play  an important
role in establishing Harnack inequalities and heat kernel
estimates  in \cite{FMT}  in the isotropic case,
whereas   in the non isotropic case,
 weighted Sobolev  inequalities are crucial in establishing
Liouville type Theorems, see \cite{CM}.

This paper is organized as follows. In Sections 2, 3  and 4 we consider
the case of codimension $k=1$ case. In particular, in Section 2 we study the  case
$\xs<2$,  in Section 3 the critical case  $\xs=2$, whereas in Section 4
the supercritical case $\xs>2$. Finally the last Section 5 is devoted
to the study of the higher  codimension case  and in particular we give the
proofs of Theorems   \ref{th9} and \ref{th40}.

 \end{section}

\begin{section}{Codimension $1$ degeneracy; the case $\xs<2$.}

In this Section we will give the proofs of
 Theorems \ref{th1}, \ref{th2} and Corollary \ref{th3}.

We first give the proof Theorem \ref{th2}.

\noindent {\em Proof of Theorem \ref{th2}:} Let us first give the
proof of part (ii). For any $u\in C^\infty_0(\R^N)$ it is well
known that
\begin{equation}\label{settebis}
S_N ||u||_{L^{\frac{N}{N-1}}} \le ||\nabla
u||_{L^1} \ ,
\end{equation}
 where
$S_N:=N\pi^{\frac{1}{2}}\left[\Gamma(1+\frac{N}{2})\right]^{-\frac{1}{2}}$
(see, e.g.,  p. 189 in \cite{M}).  We apply (\ref{settebis}) to
the function $u:=x_N^a  v$, for  $v\in
C^\infty_0(\R^N)$ and $a>0$.  Thus,  we have
$$
S_N ||x_N^a v||_{L^{\frac{N}{N-1}}} \le \int_{\R^N_{+}}
 \left(|\nabla v|x_N^a+a x_N^{a-1} |v|\right) dx' dx_N  \ .
$$
 To estimate the last term of the
right hand side, we  integrate by parts, \be\la{5}
a\int_{\R^N_{+}} x_N^{a-1} |v|dx' dx_N = \int_{\R^N_{+}} \nabla
x_N^{a} |v| dx' dx_N =- \int_{\R^N_{+}} x_N^{a}\nabla |v| dx'
dx_N. \ee From this we get
\begin{equation}\label{ottobis}
a \int_{\R^N_{+}} x_N^{a-1} |v|dx' dx_N  \leq
 \int_{\R^N_{+}} |\nabla v|x_N^a dx'
dx_N.\end{equation}
Consequently,
\begin{equation}\label{novebis}||x_N^a v||_{L^{\frac{N}{N-1}}} \le 2S_N^{-1}
\int_{\R^N_{+}} |\nabla v|x_N^a dx' dx_N .\end{equation}

For any $1\le p \le \frac{N}{N-1}$ and any two functions $w$ and
$v$,
 the following interpolation
inequality can be easily seen to be true:
\begin{equation}
\label{diecibis}||w^b v||_{L^p}\le C_1 ||w^a
v||_{L^{\frac{N}{N-1}}} + C_2||w^{a-1} v||_{L^1} \ ,
 \ \hbox { for }~~ b=a-1+\frac{p-1}{p}N,
\end{equation}
with  two positive
constants $C_1, C_2$ independent  of $w$ and $v$.

 From  (\ref{novebis}) and
(\ref{diecibis}) for  $w:=x_N$ we obtain  the following
\begin{equation}
\label{7} \left(\int_{\R^N_{+}} x_N^{bp} |v|^p dx' dx_N
\right)^\frac{1}{p} \le C_1 \int_{\R^N_{+}}|\nabla v|x_N^a dx' dx_N
+ C_2 \int_{\R^N_{+}} x_N^{a-1} |v|dx' dx_N.
\end{equation}
Using now (\ref{ottobis}) we arrive at the  following $L^p-L^1$
 weighted estimate
\begin{equation}
\label{undicibis} \left(\int_{\R^N_{+}} x_N^{bp} |v|^p dx' dx_N
\right)^\frac{1}{p} \le C_1 \int_{\R^N_{+}}|\nabla v|x_N^a dx' dx_N .
\ee
To pass to the corresponding $L^Q-L^2$ estimate
 we apply (\ref{undicibis}) to $v:=|f|^s$, $s>0$, to  obtain
$$\left(\int_{\R^N_{+}} x_N^{bp} |f|^{ps} dx' dx_N \right)^\frac{1}{p}\le
 C \int_{\R^N_{+}}
 f^{s-1}|\nabla f|x_N^{a} dx' dx_N = $$ $$=C \int_{\R^N_{+}} x_N^{\frac{bp}{2}}
 |f|^{s-1} |\nabla f| x_N^{a-\frac{bp}{2}} dx' dx_N\le $$
$$\le C \left(\int_{\R^N_{+}} x_N^{bp}
|f|^{2s-2}  dx' dx_N\right)^\frac{1}{2}\left(\int_{\R^N_{+}}
|\nabla f|^2 x_N^{2a-bp} dx' dx_N\right)^\frac{1}{2}.
$$
Choosing $s=\frac{2}{2-p}$ so that $2s-2=ps$ we get
\begin{equation}\label{dodicibis}
\left(\int_{\R^N_{+}} x_N^{bp} |f|^{ps} dx'
dx_N\right)^{\frac{2}{p}-1}\le C\int_{\R^N_{+}} x_N^{2a-bp}  |\nabla f|^2
 dx' dx_N.
\end{equation}
To  arrive at (\ref{due}) we take $BQ=bp$, $Q=ps$ and $2a-bp=2A$.
For this choice of the parameters we arrive at
$$\left(\int_{\R^N_{+}}  x_N^{BQ} |f|^{Q} dx' dx_N
\right)^\frac{2}{Q} \le C \int_{\R^N_{+}}  x_N^{2A}  |\nabla f|^2
dx' dx_N\, $$ with $2\le Q \le \frac{2N}{N-2}$ and    $B=A-1+
\frac{Q-2}{2Q}N$, in case   $N \geq 3$,  or  $Q \ge 2$ and
$B=A-\frac{2}{Q}$   in case $N=2$.

Since $2a=2A+BQ$, the condition $a>0$  is equivalent to $BQ+2A>0$.
This completes the proof of part (ii)  of   Theorem  \ref{th2}.

Concerning part (i), we note that for $v \in C^{\infty}_0(\R^N_{+})$,
and $a \in \R$,  it follows from (\ref{5}) that
\begin{equation}\label{8c}
|a| \int_{\R^N_{+}} x_N^{a-1} |v|dx' dx_N  \leq
 \int_{\R^N_{+}} |\nabla v|x_N^a dx'
dx_N.
\end{equation}
Consequently, estimate   (\ref{novebis})   remains true   for any $a \in \R$.
Estimate (\ref{7}) is still true, and using (\ref{8c}) we arrive at (\ref{undicibis}).
The use of (\ref{8c}) however imposes the condition
 that $a \neq 0$. The rest of the argument remains the same. The condition $a \neq 0$
is equivalent to $BQ+2A \neq 0$.

We  finally note that, when
$A=B=0$ then (\ref{due}) is  the standard Sobolev
inequality.

\finedim

As a consequence of the Theorem \ref{th2} we have the following
inequality in a strip:

\bigskip
\begin{proposition}\label{cor2}
 Let $\mathcal H_1   = \{(x',x_N)\in \R^{N-1}\times \R: |x'|<1\}$,
\[
N=2,  ~~~~~~~ 2  \le Q,   ~~~~~~~~{\rm and}~~~~~~ B=A-
\frac{2}{Q},
\]
or
\[ N\geq 3, ~~~~~~~ 2  \le Q  \le \frac{2N}{N-2}, ~~~~~~{\rm
and}~~~~~~ B=A-1+\frac{Q-2}{2Q}N.
\]
If  $BQ+2A \neq 0$, or if $A=B=0$ then, \\
(i)  There exists a positive constant $C=C(A,Q,N)$, such that for
any function $f\in C^\infty_0(\mathcal H_1 \cap \R^N_{+})$ there holds
\begin{equation}\label{tre}
 \left(\int_{\mathcal H_1\cap\{0<x_N<1\}} x_N^{BQ}
|f(x',x_N)|^Q dx' dx_N \right)^{\frac{2}{Q}} \le C \int_{\mathcal H_1\cap
\{0<x_N<1\}} x_N^{2A} \left(|\nabla_{x'} f|^2 +|\partial_{x_N}
f|^2\right) \ dx' dx_N \ ,
\end{equation}
(ii)
If moreover  $BQ+2A>0$,
 inequality  (\ref{tre}) still holds even if  $f  \in
C^\infty_0(\mathcal H_1)$
\end{proposition}

In the  case where $2A=BQ \in (0,\infty)$ and under the more
restrictive assumption that $f\in C^\infty_0(\mathcal H_1\cap
\{0<x_N<1\})$, the result of part (ii) has been established in [C]
by different methods (see also \cite{CF}).

\medskip

\noindent{\em Proof of Proposition  \ref{cor2}:} We prove part
(ii), the other case  being quite similar. In order to do this due
to Theorem \ref{th2} part (ii) it is enough to remove zero
boundary conditions on the hyperplane $x_N=1$. Let $f\in
C^\infty_0(\mathcal H_1)$ and we denote by $\xi(x_N)$  a $C^1$
function such that $\xi(x_N)=1$ if $x_N\le \frac{1}{2}$ and
$\xi(x_N)=0$ if $x_N\ge 1$. We then have \bea\la{2.1} LHS &  := &
C \left(\int_{\{0<x_N<1\}} x_N^{BQ} |f|^Q dx' dx_N
\right)^{\frac{2}{Q}}   \nonumber \\
 &  \le  &  \left(\int_{\{0<x_N<1\}} x_N^{BQ}
|f\xi|^Q dx' dx_N \right)^{\frac{2}{Q}} +
\left(\int_{\{\frac{1}{2}<x_N<1\}} x_N^{BQ} |f(1-\xi)|^Q dx'
dx_N \right)^{\frac{2}{Q}} \nonumber \\
& =: & I_1 + I_2.
\eea
 Applying Theorem \ref{th2} part (ii) to
the function $f\xi$, we obtain \bea\la{2.2}
 I_1 &  \le  &  C    \int_{\{0<x_N<1\}} x_N^{2A} \left(|\nabla_{x'}
(f\xi)|^2 +|\partial_{x_N} (f\xi)|^2\right) \ dx' dx_N
\nonumber \\
 &  \le  &  C
\int_{\{0<x_N<1\}} x_N^{2A} \left(|\nabla_{x'} f|^2 +|\partial_{x_N}
f|^2  + f^2\right) \ dx' dx_N. \eea Concerning $I_2$ we note that
the weights
 $x_N^{BQ}$ and  $x_N^{2A}$ are  uniformly
bounded both  from above and below for  $x_N \in [\frac{1}{2},1]$,
and therefore,
 applying
the standard Sobolev inequality to the function $f(1-\xi)$ which
is zero for  $|x'|=1$ as well as for  $x_N=\frac{1}{2}$ we get
\[
I_2 \leq C
 \int_{\{\frac{1}{2}<x_N<1\}} x_N^{2A}\left(|\nabla_{x'}
f|^2 +|\partial_{x_N} f|^2 +f^2 \right) \ dx' dx_N.
\]
Combining this with (\ref{2.1}) and (\ref{2.2}) we  get
\be\la{2.3}
LHS \leq
 C
\int_{\{0<x_N<1\}} x_N^{2A} \left(|\nabla_{x'}f|^2
 +|\partial_{x_N} f|^2  + f^2\right) \ dx' dx_N.
\ee
To continue, let
 $B_1':=\{x'\in \R^{N-1}:|x'|<1\}$. For any fixed $x_N \in [0,1]$, we have
by the
Poincar\'e inequality
$$\int_{B_1'} f^2(x',x_N) dx'\le C \int_{B_1'} |\nabla_{x'} f|^2 dx',$$
whence
$$\int^1_0\int_{B_1'} f^2(x',x_N) dx'x_N^{2A}dx_N\le C \int^1_0\int_{B_1'}
 |\nabla_{x'} f|^2 dx'x_N^{2A} dx_N \ .$$
From this and (\ref{2.3}) the result follows. \finedim

\bigskip

\noindent We are now ready to prove Theorem \ref{th1}.\\

\noindent {\em Proof of Theorem \ref{th1}:} It is enough to prove
(\ref{11}) in the upper half cylinder; that is, if $f\in
C^\infty_0(\mathcal H_1)$  then we will show that
\begin{equation} \label{rm1}
 \left(\int_{\{0<\l<1\}} (1-\lambda)^\alpha |f(x',\lambda)|^Q
dx' d\lambda \right)^{\frac{2}{Q}} \le C
\int_{\{0<\l<1\}}(1-\lambda)^\alpha \left(|\nabla_{x'} f|^2
+(1-\lambda)^{\sigma}|\partial_\lambda f|^2\right) dx' d\lambda \ .
\end{equation}
We first consider the case $\sigma< 2$. We change variables by
$x'=x' \ , s=(1-\lambda)^{\frac{2-\sigma}{2}}$ thus setting
$\varphi(x',s):=f(x', 1-s^\frac{2}{2-\sigma})$, it follows that
 inequality (\ref{rm1}) is equivalent to
\begin{equation}\label{milu}
\left(\int_{\{0<s<1\}} s^{\frac{\sigma+2\a}{2-\sigma}}
|\varphi(x',s)|^Q dx' ds \right)^{\frac{2}{Q}} \le C
\int_{\{0<s<1\}} s^{\frac{\sigma+2\a}{2-\sigma}}
\left(|\nabla_{x'} \varphi|^2 +|\partial_{s} \varphi|^2\right) \
dx' ds \ ,
\end{equation}
in fact we easily compute that $ds=
\frac{\sigma-2}{2}(1-\lambda)^{-\frac{\sigma}{2}}
  d\lambda$,
 $\partial_\lambda = \frac{ds}{d\lambda} \partial_s=
 \frac{\sigma-2}{2}(1-\lambda)^{-\frac{\sigma}{2}}
 \partial_s$ and
$$|\nabla_{x'} f|^2
+(1-\lambda)^\sigma|\partial_\lambda f|^2= |\nabla_{x'}\varphi|^2+
\left( \frac{\sigma-2}{2}\right)^2|\partial_{s} \varphi|^2 \ .$$

When $\sigma\in (-2\alpha,2)$ we now use Proposition \ref{cor2},
part (ii). Suppose first that $N \geq 3$.
 For $A=\frac{\sigma+2\a}{2(2-\sigma)}$
and  $B= \frac{\sigma+2\a}{2(2-\sigma)}-1 + \frac{Q-2}{2Q}N$, with  $2 \leq
Q \leq \frac{2N}{N-2}$ we have that the right hand side of (\ref{milu})
dominates
\[
\left(\int_{\{0<s<1\}}
s^{BQ} |\varphi(x',s)|^Q dx' ds
\right)^{\frac{2}{Q}}.
\]
To deduce  (\ref{milu})  we   need  $\frac{\sigma+2 \a}{2-\sigma}
 \geq  BQ=\left(\frac{\sigma+2\a}{2(2-\sigma)}-1 + \frac{Q-2}{2Q}N \right)Q  $,
which is equivalent to
\be\la{con1}
 Q \leq  2 \frac{ N + \frac{2 \a +\xs}{2-\xs}}{N +
\frac{2 \a +\xs}{2-\xs} -2}. \ee On the other hand, the
restriction $2A+BQ>0$ is easily seen to be equivalent to
\be\la{con2} Q >2 \frac{ N - \frac{2 \a +\xs}{2-\xs}}{N + \frac{2
\a +\xs}{2-\xs} -2}=:\bar Q. \ee We note that $Q_{cr}$ as given by
(\ref{qcr}) satisfies both (\ref{con1}) and (\ref{con2}) and
therefore  (\ref{11}) has been proved  for $Q=Q_{cr}$. The full
range of $Q$ follows by   using Holder's inequality in the left
hand side of (\ref{11}).

The case $N =2$ is treated quite similarly. Thus (\ref{11}) has
been proved  for any $f\in C^\infty_0(\mathcal H_1)$.

In the special case $\xs=2$  and $Q=2$ we note
that (\ref{11}) is still valid.
To see this we change variables by $x'=x'$ and
$t=(1-\lambda)^{\alpha+1}$ thus setting $g(x',t):=f(x',
1-t^{\frac{1}{\alpha+1}})$. It follows that
inequality (\ref{rm1}) is equivalent to
\begin{equation} \label{holi}
\left(\int_{t\in (0,1)} |g(x',t)|^Q dx' dt \right)^{\frac{2}{Q}} \le C
\int_{t\in (0,1)} \left(|\nabla_{x'} g|^2 + t^2 |\partial_{t}
g|^2\right) \ dx' dt \ ,
\end{equation}
in fact we easily compute that $dt=-(\alpha+1)(1-\lambda)^{\alpha}
  d\lambda$,
 $\partial_\lambda = \frac{dt}{d\lambda} \partial_t=  -(\alpha+1)(1-\lambda)^{\alpha}
 \partial_t$ and
$$|\nabla_{x'} f|^2
+(1-\lambda)^\sigma|\partial_\lambda f|^2= |\nabla_{x'}g|^2+
\left(\alpha+1\right)^2 t^2|\partial_{t} g|^2 \ .$$

\noindent Inequality (\ref{holi}) with $Q=2$ holds true, as one
can easily see using Poincar\'e inequality for the slices $t=$
constant.

It remains to show that (\ref{11}) fails  in the  case $\xs=2$,
$\alpha > -1$ and $Q>2$ even thought we take $f\in
C^\infty_0(\mathcal C_1)$. To this end, let us make use of the
following different change of variables $x'=x'$ and $\lambda=\tanh
x_N$. Then $\lambda  \in (-1,1)$ goes to $x_N \in (-\infty,
\infty)$ and
 $(1-|\lambda|)\sim (1-\lambda^2)=(\cosh x_N)^{-2}\sim e^{-2|x_N|}$ and
 $d\lambda  \sim  (\cosh x_N)^{-2}
dx_N \sim  e^{-2|x_N|}  dx_N  $. We define $g(x',x_N):=f(x',\tan
hx_N)$.
 Then it  follows from that for any function
 $g \in C_0^{\infty} (\mathcal H_1)$ the following inequality should be true if (\ref{11}) holds true:
\begin{equation} \label{cinque}
 \left(\int_{\mathcal H_1}
e^{-2(\alpha+1)|x_N|} |g(x',x_N)|^Q dx' dx_N \right)^{\frac{2}{Q}}
\le C \int_{\mathcal H_1}  e^{-2(\alpha+1)|x_N|}   \left(|\nabla_{x'} g|^2
+|\partial_{x_N} g|^2\right) dx' dx_N \ .
\end{equation}
For  $g \in C_0^{\infty} (\mathcal H_1 \cap \{x_N>0\})$ we set
$g_{\tau}(x',x_N):=g(x', x_N- \tau)$,   $\tau >0$. Clearly,
$g_{\tau}  \in C_0^{\infty} (\mathcal H_1 \cap \{x_N>0\})$ and
 applying (\ref{cinque}) to
the family $g_{\tau}$ we get
\[
\left(\int_{\mathcal H_1} e^{-2(\alpha+1)x_N} |g(x',x_N)|^Q dx' dx_N
\right)^{\frac{2}{Q}} \le C
e^{-2\tau(\alpha+1)\left(\frac{Q-2}{Q}\right)} \int_{\mathcal H_1}
e^{-2(\alpha+1)x_N} \left(|\nabla_{x'} g|^2 +|\partial_{x_N}
g|^2\right) dx' dx_N \ ,
\]
for any $\tau >0$. Taking the limit $\tau \to +\infty$  we reach a
contradiction for $Q>2$, $\alpha>-1$.

This completes the proof of Theorem \ref{th1}.

\finedim

\medskip

\noindent {\bf  Remark.} In case  $\sigma=2\alpha$ and $\alpha\in
(0,1)$,  estimate (\ref{11})  is an improvement of Proposition 4
of Caffarelli and Cordoba \cite{CC}. Indeed,
 our Sobolev exponent $Q_{cr}$
 is strictly  bigger than the one coming from the
arguments of \cite{CC} -- which is less than
$\frac{2N}{N+\frac{4\alpha}{\alpha+1}-2}$.
 Moreover, we only assume
that $f\in C^\infty_0(\mathcal H_1)$  instead of $f\in
C^\infty_0(\mathcal C_1)$.

\noindent {\bf Remark.} In case $\sigma=-\alpha$ and $\alpha>0$
inequality (\ref{11}) is a Sobolev inequality for a Grushin type
operator corresponding to the  vector fields
$((1-|\lambda|)^{\frac{\alpha}{2}} \nabla_{x'},
\partial_\lambda)$; we refer to \cite{FL2} where  local versions of
  similar inequalities  have been considered.

\noindent {\bf Remark.} We note that in the case
$\sigma=-2\alpha$, estimate (\ref{milu}) corresponds to the
standard Sobolev inequality in a strip, and the result follows
from Proposition \ref{cor2} part (i); thus (\ref{11}) still holds
true for any $f\in C^\infty_0(\mathcal C_1)$ if $\sigma=-2\alpha$.

\medskip

We next show  how  Corollary \ref{th3}  follows  from
Theorem \ref{th1}.

\bigskip
\noindent {\em Proof of Corollary \ref{th3}:} It is a consequence of
Theorem \ref{th1}. Indeed, for $(x', \lambda) \in C_{1,\epsilon}$
 we have $1-|\lambda|>\epsilon^{1+m}$, that is,
$\epsilon^{-1}>(1-|\lambda|)^{-\frac{1}{1+m}}$, and so
$\frac{(1-|\lambda|)^\beta}{\epsilon^\beta}>(1-|\lambda|)^{\frac{\beta
m}{1+m}}$,  $\beta >0$. The result then  follows from Theorem \ref{th1} by
choosing $\sigma:=\frac{\beta m}{1+m}$ there; in particular $P_{cr}(N,m,
\alpha, \beta) = Q_{cr}(N,\alpha,\frac{\beta m}{1+m})$.

\finedim

\end{section}


\setcounter{equation}{0}

\begin{section}{\bf The critical case $\sigma=2$.}


As we have seen in  Theorem \ref{th1} inequality (\ref{11}) fails
for $\xs=2$,  $\alpha>-1$  and $Q>2$.  To obtain Sobolev type
inequalities in this case, we need to use different weights in the two
sides of the inequality. More
precisely we have the following

\begin{theorem} \label{prop4} Let $N\ge 2$,  and  $\alpha > -1$. For any $Q$ with
 $2 \leq  Q \le \frac{2N}{N-2}$, in case $N \geq 3$, or $Q \geq 2$ in case
 $N=2$, and for any $\theta>\frac{(Q-2)(\alpha+1)}{2}$
there exists a positive constant  $C=C(N,\alpha,Q,\theta)$, such
that
  for any function $f\in C^\infty_0(\mathcal
C_1)$  there  holds
\begin{equation}\label{1111}
\left(\int_{\mathcal C_1} (1-|\lambda|)^{\alpha+\theta}
 |f(x',\lambda)|^Q dx' d\lambda
\right)^{\frac{2}{Q}} \le C  \int_{\mathcal C_1}
(1-|\lambda|)^\alpha   \left(|\nabla_{x'} f|^2
+(1-|\lambda|)^{2}|\partial_\lambda f|^2\right) dx' d\lambda .
\end{equation}
\end{theorem}

\noindent {\em Proof :} It is enough to prove (\ref{1111}) in the
upper half cylinder; that is, if $f\in C^\infty_0(\mathcal
C_1)$ then we will show that
\begin{equation} \label{1111bis}
 \left(\int_{\{0<\l<1\}} (1-\lambda)^{\alpha+\theta} |f(x',\lambda)|^Q
dx' d\lambda \right)^{\frac{2}{Q}} \le C
\int_{\{0<\l<1\}}(1-\lambda)^\alpha \left(|\nabla_{x'} f|^2
+(1-\lambda)^{2}|\partial_\lambda f|^2\right) dx' d\lambda \ .
\end{equation}
We change variables by $x'=x'$, $s=-\frac{1}{K}\ln(1-\lambda)$,
for an arbitrary $K> 0$, thus setting
$\varphi(x',s)=f(x',1-e^{-Ks})$, and arguing as in the proof of
Theorem \ref{th1}, we see that inequality (\ref{1111bis}) follows
as soon as we prove the following inequality
\begin{equation}\label{mare}
\left(\int_{\{s>0\}} e^{-sK(\alpha+\theta+1)}
|\varphi(x',s)|^Q
 dx' ds \right)^{\frac{2}{Q}} \le C
\int_{\{s>0\}} |\nabla \varphi|^2 e^{-sK(\alpha+1)} dx' ds \ .
\end{equation}
 In fact we easily compute that $d\lambda= K
e^{-Ks} ds =K (1-\lambda) ds$,
 $\partial_\lambda =
\frac{ds}{d\lambda}\pa_s  = \frac{1}{K}(1-\lambda)^{-1}
\partial_{s} $ and $$|\nabla_{x'} f|^2
+(1-\lambda)^{2}|\partial_\lambda f|^2=|\nabla_{x'}\varphi|^2 +
\frac{1}{K^2} |\partial_{s} \varphi|^2\sim |\nabla \varphi|^2 \
.$$ We note that $\varphi\in C^\infty_0(\mathcal H_1)$.

To continue, we will make use of Proposition \ref{corlemnew} see
below.
 For  $A=\frac{K(\alpha+1)}{2}$ and
$B= \frac{K(\alpha+1)}{2} +\frac{Q-2}{2Q}N = \frac{1}{Q} \left(
K(\alpha+1)+ \frac{(N+K(\alpha+1))(Q-2)}{2} \right)$ we have
\begin{equation}\label{pas}
\left(\int_{\R^N_{+}}e^{-s K(\alpha+1)} e^{-\theta K s}
| \varphi(x',s)|^Q
 dx' ds \right)^{\frac{2}{Q}} \le C
\int_{\R^N_{+}}   e^{-s K (\alpha+1)}  |\nabla  \varphi|^2 dx' ds \
,~~~~~\forall  \varphi\in C^\infty_0(\mathcal H_1),
\end{equation}
where $\theta:=(\frac{N}{K}+\alpha+1)\frac{Q-2}{2}$.  Note that
$\theta=0$ if $Q=2$ as suggested by Theorem \ref{th1}. Due to the
arbitrariness of $K$ this means that we may take any value
$\theta>\frac{(\alpha+1)(Q-2)}{2}$. The restriction $2A+BQ \neq 0$
is easily seen to  be equivalent to $\frac{Q+2}{2} \left(
K(\alpha+1) + \frac{Q-2}{Q+2}N \right) \neq 0$, which is trivially
satisfied.

The case $N=2$ is treated quite similarly.

\finedim

\medskip

According to Theorem \ref{prop4}  one cannot match the weights
in the weighted anisotropic Sobolev inequality (\ref{1111})
when  $\alpha>-1$  and $Q>2$.  However,  in the case $\alpha<-1$
we can match the weights, thus proving Theorem \ref{th4444}.

\medskip

\noindent {\em Proof of Theorem \ref{th4444} :} The case $Q=2$ is a simple
 consequence
of Poincar\'e inequality. We therefore consider the case  $Q>2$.
Using the same change of variables as in the proof of Theorem
\ref{prop4} the sought for inequality is equivalent to the
following inequality
\begin{equation}\label{pas bis}
\left(\int_{\R^N_{+}}e^{-s K(\alpha+1)}
| \varphi(x',s)|^Q
 dx' ds \right)^{\frac{2}{Q}} \le C
\int_{\R^N_{+}}   e^{-s K (\alpha+1)}  |\nabla  \varphi|^2 dx' ds \
,~~~~~\forall  \varphi\in C^\infty_0(\mathcal H_1).
\end{equation}
We will use Proposition \ref{pr2}.
Thus, we have
\begin{equation}\label{3.1}
\left(\int_{\R^N_{+}}  e^{BQ s} | \varphi|^{Q} dx' ds
\right)^\frac{2}{Q} \le C \int_{\R^N_{+}} e^{2As}  |\nabla  \varphi|^2
 dx' ds,
\end{equation}
for $B=-\frac{K(\xa+1)}{Q}$ and $A-1=B-\frac{Q-2}{2Q}N
=-\frac{K(\xa+1)}{Q}-\frac{Q-2}{2Q}N$. To  deduce (\ref{pas bis})
from (\ref{3.1}) we need $2A \leq -K(\xa+1)$ which is equivalent
to  $2 \leq \frac{Q-2}{Q} \left(N - K (\xa+1) \right)$. This last
inequality is always satisfied by taking $K$ large enough. On the
other hand $BQ+2(A-1)= -K (\xa+1) \frac{Q+2}{Q} - \frac{Q-2}{Q}N
\neq 0$, for $K$ large.

The case $N=2$ is treated  similarly.

\finedim

It remains to give the proof  of the auxiliary results we used above.
We first have

\begin{theorem}\la{lemnew}
 Let either
\[
N=2,  ~~~~~~~ 2  \le Q,   ~~~~~~~~{\rm and}~~~~~~ B=A +1-
\frac{2}{Q},
\]
or else,
\[
N\geq 3, ~~~~~~~ 2  \le Q  \le \frac{2N}{N-2}, ~~~~~{\rm and}~~~~~
B=A +\frac{Q-2}{2Q}N.
\]
Then, if  $BQ+2A \neq 0$, there exists a positive constant $C=C(A,Q,N)$ such that for
any function  $f\in C^\infty_0(\R^N_{+})$ there holds
\begin{equation}\label{2.10}
\left(\int_{\R^N_{+}}  e^{-BQ x_N} |f|^{Q} dx' dx_N
\right)^\frac{2}{Q} \le C \int_{\R^N_{+}} e^{-2Ax_N}  |\nabla f|^2
 dx' dx_N .
\end{equation}
\end{theorem}

\medskip
\noindent {\em Proof:} We apply  the Gagliardo--Nirenberg--Sobolev
inequality (\ref{settebis}) to the function $u:=e^{-a x_N} v$, for
any $v \in C^\infty_0(\R^N_{+})$ and $a \neq 0$, to get
$$
S_N ||e^{-a x_N} v||_{L^{\frac{N}{N-1}}} \le \int_{\R^N_{+}}
 \left(|\nabla v|e^{-ax_N} +|a| e^{-ax_N}
|v|\right) dx' dx_N  \ .
$$
 To estimate the last term of the right
hand side,  we integrate  by parts,
\[
a\int_{\R^N_{+}} e^{-ax_N} |v|dx' dx_N = -\int_{\R^N_{+}} \nabla
e^{- a x_N} |v| dx' dx_N =\int_{\R^N_{+}} e^{-a x_N}\nabla |v| dx'
dx_N
\]
whence,
\begin{equation}\label{otto}
|a|\int_{\R^N_{+}} e^{-ax_N} |v|dx' dx_N  \le
\int_{\R^N_{+}}e^{-ax_N} |\nabla v| dx' dx_N.
\end{equation}
Consequently,
\begin{equation}\label{nove}||e^{-a x_N} v||_{L^{\frac{N}{N-1}}} \le C
\int_{\R^N_{+}}e^{-ax_N} |\nabla v| dx' dx_N ;
\end{equation}
We note that this is true even if $a=0$.

Using the interpolation inequality (\ref{diecibis}) with
$w:=e^{-x_N}$, as well as
 (\ref{otto}) and (\ref{nove}) we arrive at the following $L^p-L^1$ estimate
 ($e^{-(a-1)x_N}\ge e^{-ax_N}$)
\begin{equation}
\label{undici} \left(\int_{\R^N_{+}} e^{-bpx_N} |v|^p dx' dx_N
\right)^\frac{1}{p} \le C \int_{\R^N_{+}}e^{-(a-1)x_N}  |\nabla v|
dx' dx_N \ ,
\end{equation}
with  $1 \le p \le \frac{N}{N-1}$, $b=a-1 +\frac{p-1}{p}N$ and $a
\neq 1$. Indeed in order to reach inequality (\ref{undici}) we
need the following inequality
$$
\int_{\R^N_{+}} e^{-(a-1)x_N} |v|dx' dx_N  \le C
\int_{\R^N_{+}}e^{-(a-1)x_N} |\nabla v| dx' dx_N.
$$
which follows from inequality (\ref{otto}) if $a\neq 1$.

We next  apply (\ref{undici}) to $v:=|f|^s$,  $s
>0$, to  obtain
$$\left(\int_{\R^N_{+}} e^{-bpx_N} |f|^{ps} dx' dx_N \right)^\frac{1}{p}\le C \int_{\R^N_{+}}
 f^{s-1}|\nabla f|e^{-(a-1)x_N} dx' dx_N = $$ $$=C \int_{\R^N_{+}} e^{-\frac{bpx_N}{2}}
f^{s-1} |\nabla f| e^{-(a-1)x_N+\frac{bpx_N}{2}} dx' dx_N\le
$$
$$\le C \left(\int_{\R^N_{+}} e^{-bpx_N}
f^{2s-2}  d x' dx_N\right)^\frac{1}{2}\left(\int_{\R^N_{+}}
|\nabla f|^2 e^{-2(a-1)x_N+bpx_N} dx' dx_N\right)^\frac{1}{2}.
$$
Choosing $s = \frac{2}{2-p}$, so that $2s-2=ps$ we get
\begin{equation}\label{dodici}
\left(\int_{\R^N_{+}} e^{-bpx_N} |f|^{ps} dx'
dx_N\right)^{\frac{2}{p}-1}\le C\int_{\R^N_{+}}
e^{(-2(a-1)+bp)x_N} |\nabla f|^2
 dx' dx_N.
\end{equation}
To conclude the proof of the Lemma we take $BQ=bp$, $Q=ps$,
 and $A=a -1 - \frac{bp}{2}$. The condition $a\neq 1$ is equivalent
to  $BQ+2A \neq 0$.

\finedim

As a consequence of the previous Theorem, we have the following
result which is the analogue of Proposition  \ref{cor2}.
That is, is some cases we can remove the zero boundary condition
at $x_N=0$.

\begin{proposition}\la{corlemnew}
 Let either
\[
N=2,  ~~~~~~~ 2  \le Q,   ~~~~~~~~{\rm and}~~~~~~ B=A +1-
\frac{2}{Q},
\]
or else,
\[
N\geq 3, ~~~~~~~ 2  \le Q  \le \frac{2N}{N-2}, ~~~~~{\rm and}~~~~~
B=A +\frac{Q-2}{2Q}N.
\]

Then, if  $BQ+2A \neq 0$,  there exists a positive constant $C=C(A,Q,N)$
 such that for
any function  $f\in C^\infty_0(\mathcal H_1)$ there holds
\begin{equation}\label{22.10}
\left(\int_{\R^N_{+}}  e^{-BQ x_N} |f|^{Q} dx' dx_N
\right)^\frac{2}{Q} \le C \int_{\R^N_{+}} e^{-2Ax_N}  |\nabla f|^2
 dx' dx_N .
\end{equation}
\end{proposition}

\medskip
\noindent {\em Proof. } To deduce (\ref{22.10}) from
(\ref{2.10}) we will work as in the proof of Proposition  \ref{cor2}
in order to remove the zero boundary condition on the hyperplane
$x_N=0$.
 Let $\xi(x_N)$ be a $C^1$ function such that $\xi(x_N)=1$ if
$x_N \geq
 2$ and $\xi(x_N)=0$ if $x_N\in[0,1]$, then for any $f\in
C^\infty_0(\mathcal H_1)$ we have \bea LHS & := & C
\left(\int_{\R^N_{+}}e^{-BQx_N} |f|^Q dx'
dx_N \right)^{\frac{2}{Q}}  \nonumber \\
&  \le  &   \left(\int_{\R^N_{+}} e^{-BQx_N} |f\xi|^Q dx' dx_N
\right)^{\frac{2}{Q}}+
 \left(\int_{\R^N_{+}}
e^{-BQx_N} |f(1-\xi)|^Q   dx' dx_N \right)^{\frac{2}{Q}} \nonumber  \\
& =: & I_1 + I_2. \eea Applying (\ref{2.10}) to the function
$f\xi$, we obtain
\[
I_1 \le C \int_{\R^N_{+}} e^{-2Ax_N} \left(|\nabla f|^2 + f^2
\right) \ dx' dx_N.
\]
On the other hand, since  the weights $e^{-BQx_N}$ and
$e^{-2Ax_N}$ are uniformly bounded both from above and below in
the interval $[0,2]$, we may apply the
 standard Sobolev inequality to the function $f(1-\xi)$ which
is zero when $|x'|=1$ as well as when $x_N=2$ to  get
\[
I_2 \leq C
 \int_{\R^{N-1}\times[0,2]} |\nabla
(f(1-\xi))|^2  \ dx' dx_N  \leq
 C \int_{\R^N_{+}} e^{-2Ax_N}
\left(|\nabla f|^2 + f^2 \right) \ dx' dx_N.
\]
Combining the above estimates we have \be\la{2.15bis} LHS \leq C
\int_{\R^N_{+}} e^{-2Ax_N} \left(|\nabla f|^2 + f^2 \right) \ dx'
dx_N. \ee To conclude we  use the Poincar\'e inequality
 on the set $B_1'=\{x'\in
\R^{N-1}: |x'|<1\}$. For any fixed $x_N$
$$\int_{B_1'} f^2(x',x_N) dx'\le C \int_{B_1'} |\nabla_{x'} f|^2 dx',$$
whence
$$
\int_{\R^N_{+}} e^{-2Ax_N} f^2 dx' dx_N = \int_0^\infty e^{-2Ax_N}
\int_{B_1'} f^2 dx' dx_N\le C \int_0^\infty e^{-2A x_N}
\int_{B_1'} |\nabla_{x'} f|^2 dx' dx_N \ ,$$ $$\leq C
\int_{\R^N_{+}} e^{-2Ax_N} |\nabla f|^2 dx' dx_N .$$ From this and
(\ref{2.15bis}) the result follows.

\finedim

We next present a  new Sobolev inequality which also involves exponential
weights. We used this estimate in the proof of Theorem \ref{th4444}.

\begin{theorem}\la{lemnewbis}
 Let either
\[
N=2,  ~~~~~~~ 2  \le Q,   ~~~~~~~~{\rm and}~~~~~~ B=A-
\frac{2}{Q},
\]
or else,
\[
N\geq 3, ~~~~~~~ 2  \le Q  \le \frac{2N}{N-2}, ~~~~~{\rm and}~~~~~
B=A -1+\frac{Q-2}{2Q}N.
\]
Then, if  $BQ+2A \neq 2$, there exists a positive constant $C=C(A,Q,N)$ such that for
any function  $f\in C^\infty_0(\R^N_{+})$ there holds
\begin{equation}\label{2.10bis}
\left(\int_{\R^N_{+}}  e^{BQ x_N} |f|^{Q} dx' dx_N
\right)^\frac{2}{Q} \le C \int_{\R^N_{+}} e^{2Ax_N}  |\nabla f|^2
 dx' dx_N .
\end{equation}
\end{theorem}

\medskip
\noindent {\em Proof. }
Working as in the proof of Theorem \ref{lemnew} we obtain (\ref{otto}) and
 (\ref{nove})
that is,
\begin{equation}\label{8}
|a|\int_{\R^N_{+}} e^{ax_N} |v|dx' dx_N  \le
\int_{\R^N_{+}}e^{ax_N} |\nabla v| dx' dx_N,
\end{equation}
and
\begin{equation}\label{9}||e^{a x_N} v||_{L^{\frac{N}{N-1}}} \le C
\int_{\R^N_{+}}e^{ax_N} |\nabla v| dx' dx_N ;
\end{equation}
which are valid for any $a$ in $\R$.

We next use  the interpolation inequality (\ref{diecibis}) with
$w:=e^{x_N}$, as well as
 (\ref{8}) and (\ref{9}) to  arrive at the following $L^p-L^1$ estimate
 ($e^{ax_N}\ge e^{(a-1)x_N}$)
\begin{equation}
\label{11b} \left(\int_{\R^N_{+}} e^{bpx_N} |v|^p dx' dx_N
\right)^\frac{1}{p} \le C \int_{\R^N_{+}}e^{a x_N}  |\nabla v| dx'
dx_N \ ,
\end{equation}
with  $1 \le p \le \frac{N}{N-1}$, $b=a-1 +\frac{p-1}{p}N$ and $a
\neq 1$. To  reach inequality (\ref{11b}) we  used
the following estimate
$$
\int_{\R^N_{+}} e^{(a-1)x_N} |v|dx' dx_N  \le C
\int_{\R^N_{+}}e^{(a-1)x_N} |\nabla v| dx' dx_N,
$$
which is a consequence of  (\ref{8}) if $a\neq 1$.

We next  apply (\ref{11b}) to $v:=|f|^s$,  $s
>0$, to  obtain
$$\left(\int_{\R^N_{+}} e^{bpx_N} |f|^{ps} dx' dx_N \right)^\frac{1}{p}\le C
 \int_{\R^N_{+}}
 f^{s-1}|\nabla f|e^{a x_N} dx' dx_N = $$ $$=C \int_{\R^N_{+}}
 e^{\frac{bpx_N}{2}}
f^{s-1} |\nabla f| e^{a x_N-\frac{bpx_N}{2}} dx' dx_N\le
$$
$$\le C \left(\int_{\R^N_{+}} e^{bpx_N}
f^{2s-2}  d x' dx_N\right)^\frac{1}{2}\left(\int_{\R^N_{+}}
|\nabla f|^2 e^{2a x_N-bpx_N} dx' dx_N\right)^\frac{1}{2}.
$$
Choosing $s = \frac{2}{2-p}$, so that $2s-2=ps$ we get
\begin{equation}\label{12}
\left(\int_{\R^N_{+}} e^{bpx_N} |f|^{ps} dx'
dx_N\right)^{\frac{2}{p}-1}\le C\int_{\R^N_{+}} e^{(2a-bp)x_N}
|\nabla f|^2
 dx' dx_N.
\end{equation}
To conclude the proof of the Lemma we take $BQ=bp$, $Q=ps$,
 and $A=a - \frac{bp}{2}$. The condition $a \neq 1$
is equivalent to $BQ+2A \neq 2$.

\finedim

We finally have

\begin{proposition}\la{pr2}
 Let either
\[
N=2,  ~~~~~~~ 2  \le Q,   ~~~~~~~~{\rm and}~~~~~~ B=A -
\frac{2}{Q},
\]
or else,
\[
N\geq 3, ~~~~~~~ 2  \le Q  \le \frac{2N}{N-2}, ~~~~~{\rm and}~~~~~
B=A-1 +\frac{Q-2}{2Q}N.
\]

Then, if  $BQ+2A \neq 2$,  there exists a positive constant $C=C(A,Q,N)$
 such that for
any function  $f\in C^\infty_0(\mathcal H_1)$ there holds
\begin{equation}\label{2.10bisb}
\left(\int_{\R^N_{+}}  e^{BQ x_N} |f|^{Q} dx' dx_N
\right)^\frac{2}{Q} \le C \int_{\R^N_{+}} e^{2Ax_N}  |\nabla f|^2
 dx' dx_N .
\end{equation}
\end{proposition}

{Proof:} We need    to remove the zero boundary condition of $f$,
 on the hyperplane
$x_N=0$. As usual,
 let $\xi(x_N)$ be a $C^1$ function such that $\xi(x_N)=1$ if
$x_N \geq
 2$ and $\xi(x_N)=0$ if $x_N\in[0,1]$, then for any $f\in
C^\infty_0(\mathcal H_1)$ we have  $f = f \xi+ f(1-\xi)$.
To conclude the proof we argue as in the proof of Proposition
\ref{corlemnew}.  We omit further details.

\end{section}

\setcounter{equation}{0}

\begin{section}{\bf  The supercritical  case  $\xs >2$}

In this Section we will  give the  proof of Theorem \ref{th1.3}.
It is a direct consequence of a  more
general result. We recall that
\[
Q_{cr} =
 Q_{cr}(N,\alpha,\sigma):=
\frac{2 \left( N + \frac{2 \a +\xs}{2-\xs}\right)}{N + \frac{2 \a
+\xs}{2-\xs} -2}. \]
We also set
\[
 \bar{Q}(N,\alpha,\sigma):=
\frac{2 \left( N - \frac{2 \a +\xs}{2-\xs}\right)}{N + \frac{2 \a
+\xs}{2-\xs} -2}, \]
and
\be\la{defth}
\theta_{cr}:=\frac{2-\sigma}{2}\left[ \frac{Q}{2}\left\{N+\frac{2\alpha+\sigma}{2-\sigma}
-2\right\}
 -\left\{N+\frac{2\alpha+\sigma}{2-\sigma} \right\} \right]
= \frac{2-\xs}{4}
 \left( N + \frac{2\alpha +\sigma}{2-\sigma} -2 \right) (Q - Q_{cr}) .
\ee
We then have

\begin{theorem}\label{th4.1}
 Let $N\ge 2$,
$\alpha<-1$  and   $\sigma \in (2,-2\alpha)$.
Then, for any  $\theta \geq \theta_{cr}$ and
 any $Q \neq \bar{Q}$
with $2 \leq  Q \leq \frac{2 N}{N  -2}$, in case
 $N\ge 3$, or $Q\ge 2$ in case $N=2$,
 there exists a positive constant $C=C(Q,N,\alpha,\sigma, \theta)$,
such that
 for any function $f\in C^\infty_0(\mathcal C_1)$ there holds
\begin{equation}\la{44.1}
 \left(\int_{\mathcal C_1} (1-|\lambda|)^{\alpha+\theta}
|f(x',\lambda)|^Q dx' d\lambda \right)^{\frac{2}{Q}} \le C
\int_{\mathcal C_1}  (1-|\lambda|)^\alpha \left(|\nabla_{x'} f|^2
+(1-|\lambda|)^{\sigma}|\partial_\lambda f|^2\right) dx' d\lambda
\ .
 \end{equation}
\end{theorem}

To prove the above result we will use the following consequence
of Theorem \ref{th2}

\begin{proposition}\label{cor20}
 Let $\mathcal H_1   = \{(x',x_N)\in \R^{N-1}\times \R: |x'|<1\}$,
\[
N\geq 3, ~~~~~~~ 2  \le Q  \le \frac{2N}{N-2}, ~~~~~~{\rm
and}~~~~~~ B=A-1+\frac{Q-2}{2Q}N,
\]
or
\[
N=2,  ~~~~~~~ 2  \le Q,   ~~~~~~~~{\rm and}~~~~~~ B=A- \frac{2}{Q}.
\]
If  $BQ+2A \neq 0$, or if $A=B=0$ then,
 there exists a positive constant $C=C(A,Q,N)$, such that for
any function $f\in C^\infty_0(\mathcal H_1 )$ there holds
\begin{equation}\label{44.3}
 \left(\int_{\{x_N >1\}} x_N^{BQ}
|f(x',x_N)|^Q dx' dx_N \right)^{\frac{2}{Q}} \le C \int_{
\{x_N> 1\}} x_N^{2A} \left(|\nabla_{x'} f|^2 +|\partial_{x_N}
f|^2\right) \ dx' dx_N \ .
\end{equation}
\end{proposition}

\noindent {\em Proof of Proposition  \ref{cor20}:} The proof is
quite similar to the proof of  Proposition  \ref{cor2}  we
therefore sketch it. We use a $C^1$  cutoff  function $\xi(x_N)$
such that $\xi(x_N)=1$ in $x_N\ge 2$ and $\xi(x_N)=0$ if $0\le
x_N\le 1$. Hence we write $f = f \xi+ f(1-\xi)$. Now $f(1-\xi)$
satisfies the standard Sobolev inequality in $1\le x_N\le 2$,
while $f\xi$ satisfies the assumptions  of  Theorem \ref{th2} part
(i). Putting things together and using Poincar\'e inequality in
the $x'$--variables we conclude the proof. We omit further
details.

\finedim

\noindent {\em Proof of Theorems \ref{th4.1}  and  \ref{th1.3}:}
We first prove Theorem  \ref{th4.1}. As usual,
it is enough to prove (\ref{44.1}) in the upper half cylinder. That
is, if $f\in C^\infty_0(\mathcal C_1)$  then we
need to  show that
\begin{equation} \label{rm1bis}
 \left(\int_{\{0<\l<1\}} (1-\lambda)^{\alpha+\theta} |f(x',\lambda)|^Q
dx' d\lambda \right)^{\frac{2}{Q}} \le C
\int_{\{0<\l<1\}}(1-\lambda)^\alpha \left(|\nabla_{x'} f|^2
+(1-\lambda)^{\sigma}|\partial_\lambda f|^2\right) dx' d\lambda \
.
\end{equation}
As in the proof of Theorem \ref{th1}, we change variables by
$x'=x' \ , s=(1-\lambda)^{\frac{2-\sigma}{2}}$ thus setting
$\varphi(x',s):=f(x', 1-s^\frac{2}{2-\sigma})$, it follows that
 inequality (\ref{rm1bis}) is equivalent to
\begin{equation}\label{milubis}
\left(\int_{\{s>1\}} s^{\frac{\sigma+2\a+2\theta}{2-\sigma}}
|\varphi(x',s)|^Q dx' ds \right)^{\frac{2}{Q}} \le C
\int_{\{s>1\}} s^{\frac{\sigma+2\a}{2-\sigma}} \left(|\nabla_{x'}
\varphi|^2 +|\partial_{s} \varphi|^2\right) \ dx' ds \ ,
\end{equation}
for  $ \varphi \in C^\infty_0( \mathcal H_1)$. We now use
Proposition \ref{cor20}. Suppose first that $N \geq 3$.
 For $A=\frac{\sigma+2\a}{2(2-\sigma)}$
and  $B= \frac{\sigma+2\a}{2(2-\sigma)}-1 + \frac{Q-2}{2Q}N$, with
$2 \leq Q \leq \frac{2N}{N-2}$ we have that the right hand side of
(\ref{milubis}) dominates
\[
\left(\int_{\{s>1\}} s^{BQ} |\varphi(x',s)|^Q dx' ds
\right)^{\frac{2}{Q}}.
\]
To deduce  (\ref{milubis})  we   need  $\frac{\sigma+2
\a+2\theta}{2-\sigma}
 \leq  BQ=\left(\frac{\sigma+2\a}{2(2-\sigma)}-1 + \frac{Q-2}{2Q}N \right)Q  $,
which is  satisfied by  any  $\theta  \geq \theta_{cr}$  as defined  in (\ref{defth}).

Let us finally observe that $BQ+2A   \neq 0$ corresponds to the
assumption $\theta  \neq -\sigma-2\alpha$ that is $Q\neq \bar Q$.

The case $N=2$ is treated quite similarly.

To prove Theorem \ref{th1.3} we note that for $Q \geq Q_{cr}(N,\xa,\xs)$
we have that  $\theta_{cr} \leq 0$ and therefore
we  can  take $\theta=0$.

\finedim

\end{section}

\setcounter{equation}{0}

\begin{section}{\bf The case of codimension $k$ degeneracy
 $1 <   k  < N$.}


In this section we will  prove Theorems \ref{th9} and \ref{th40}.

\noindent
 {\em  Proof of Theorem \ref{th40}:} We will divide the proof into three
steps.

\smallskip

\noindent {\bf step  1} (The  critical $L^1$  weighted anisotropic
inequality). Suppose that either $\xb > 0$ and $u \in
C_0^{\infty}(\R^N)$  or  else   $\xb \in \R$ and  $u \in
C_0^{\infty}(\R^N_+)$. Then, for a constant $C$ depending only on
$N$ there holds: \be\la{4.1}  \left(
 \int_{\R^N_{+}} x_N^{\frac{\xb N + \xg (k-1)}{N-1}}|u|^{\frac{N}{N-1}} dx \right)^{\frac{N-1}{N}}
\leq    C(N) \int_{\R^N_{+}} x_N^{\xb} \left( |\nabla_{x',x_N}u| +
x_N^{\xg}|\nabla_{y}u|\right) dx. \ee The  proof follows closely
the standard proof of the $L^1$ Gagliardo-- Nirenberg--Sobolev
inequality.  Suppose that $\xb \neq 0$ and $u
\in C_0^{\infty}(\R^N_+)$. Let us write $x'=(x'_1, \ldots,
x'_{N-k})$ and $y=(y_1, \ldots,y_{k-1})$. We then have that for
$i=1, \ldots, N-k$,
\[
u(x)= - \int^{\infty}_{x'_i}
u_{x'_i}(x'_1,...t_i,...,x'_{N-k},x_N,y) dt_i,
\]
From which it follows easily that \be\la{4.3} x_N^{\xb}|u(x)| \leq
\int_{\R}x_N^{\xb}|u_{x'_i}| dt_i.  \ee We similarly have that
\be\la{4.5}x_N^{\xg+\xb}|u(x)| \leq
\int_{\R}x_N^{\xg+\xb}|u_{y_i}| ds_i, \ee where integration is
performed in the $y_i$--variable, $i=1,\ldots, k-1$. A  similar
argument  shows that
\[  x_N^{\xb}|u(x)| \leq
\int_{0}^{\infty}(\xi^{\xb}|u_{x_N}| + |\xb|\xi^{\xb-1}|u|)d \xi,
\]
from which it follows easily that \be\la{4.7}
 x_N^{\xb}|u(x)| \leq  2
\int_{0}^{\infty} \xi^{\xb}|u_{x_N}| d \xi \ee which is true also
if $\beta=0$.

Multiplying (\ref{4.3}), (\ref{4.5}), (\ref{4.7}) and raising to
the power $\frac{1}{N-1}$ we get \be\la{4.9} x_N^{\frac{\xb N +
\xg(k-1)}{N-1}} |u(x)|^{\frac{N}{N-1}}  \leq 2\left(
\prod_{i=1}^{N-k} \int_{\R}x_N^{\xb}|u_{x'_i}| dt_i
\left(\int_{0}^{\infty} \xi^{\xb}|u_{x_N}| d \xi \right)
\prod_{j=1}^{k-1} \int_{\R}x_N^{\xg+\xb}|u_{y_j}| ds_j
\right)^{\frac{1}{N-1}}. \ee We next
 integrate with respect to $x'_1$ and apply Holder's inequality in
 the right   hand side, then we integrate with respect to the
 $x'_2$ variable and so on until we integrate with respect to all
 variables. This  way  we  reach the following   estimate
\be\la{4.11} \int_{\R^N_+}x_N^{\frac{\xb N + \xg(k-1)}{N-1}}
|u(x)|^{\frac{N}{N-1}} dx  \leq 2 \left(  \prod_{i=1}^{N-k}
\int_{\R^N_+}x_N^{\xb}|u_{x'_i}| dx \left(\int_{\R^N_+}
x_N^{\xb}|u_{x_N}| dx \right) \prod_{j=1}^{k-1}
\int_{\R^N_+}x_N^{\xg+\xb}|u_{y_j}| dx \right)^{\frac{1}{N-1}}.
\ee To continue we use in the right hand
 side of (\ref{4.11}) the well known inequality
 \[
 \prod_{i=1}^N a_i  \leq  \frac{1}{N^N}
 \left( \sum_{i=1}^N a_i\right)^{N}, \hspace{14mm} a_i \geq 0.
 \]
We then conclude that \be\la{4.13}
 \int_{\R^N_+}x_N^{\frac{\xb N + \xg(k-1)}{N-1}}
|u(x)|^{\frac{N}{N-1}} dx  \leq  C(N) \left(
 \int_{\R^N_+} x_N^{\xb} \left( |\nabla_{x',x_N}u| +
x_N^{\xg}|\nabla_{y}u|\right) dx \right)^{\frac{N}{N-1}}, \ee
which is the sought for estimate (\ref{4.1}).

\smallskip

\noindent {\bf step 2} (The $L^p$--$L^1$  estimate). For $1\leq p
\leq \frac{N}{N-1}$ we will use the interpolation inequality
(\ref{diecibis}) with weight $$ w=x_N^{\frac{N+\gamma (k-1)}{N}},
$$
and
$$
a := \frac{\xb N+ \xg(k-1)}{N + \xg(k-1)}, \hspace{14mm}
\hspace{14mm} b=a-1+\frac{p-1}{p}N. $$ For these choices we have
that \be\la{4.14} \|w^{b} u \|_{L^p} \leq C_1 \left(
\int_{\R^N_+}x_N^{\frac{\xb N + \xg(k-1)}{N-1}}
|u(x)|^{\frac{N}{N-1}} dx\right)^{\frac{N-1}{N}} + C_2
\int_{\R^N_+}x_N^{\xb-1} |u| dx. \ee We will also make use of the
estimate \be\la{4.15} |\xb| \int_{\R^N_+} x^{\xb-1}_N |u| dx \leq
\int_{\R^N_+} x^{\xb}_N |u_{x_N}| dx  \leq  \int_{\R^N_+}
x^{\xb}_N |\nabla_{x',x_N}u| dx \ee which follows easily using an
integration by parts if $\beta \neq 0$. From (\ref{4.13}),
(\ref{4.14}), (\ref{4.15})  and using the  specific values of the
weight and the parameters we get \be\la{4.17} \|x_N^{\tilde{b}} u \|_{L^p}
\leq C \int_{\R^N_+} x_N^{\xb} \left( |\nabla_{x',x_N}u| +
x_N^{\xg}|\nabla_{y}u|\right) dx,  \ee with $\xg \in \R$,
 $\xb \neq 0$ and \be\la{4.19}
 {\tilde{b}}=\xb-1+\frac{p-1}{p}(N +\xg (k-1)). \ee

\smallskip

\noindent {\bf step  3} (The $L^Q$--$L^2$ estimate). Here we will
apply estimate (\ref{4.17}) to the function $u(x)=|f(x)|^s$ with
$s>0$. After some elementary calculations and use of Holder's
inequality we find that \be\la{4.21} \left( \int_{\R^N_+} x_N^{{\tilde{b}}p}
|f|^{sp} dx \right)^{\frac{1}{p}} \leq C
 \left( \int_{\R^N_+} x_N^{{\tilde{b}}p} |f|^{2s-2} dx \right)^{\frac12}
 \left( \int_{\R^N_+} x_N^{2 \xb -{\tilde{b}}p}  \left( |\nabla_{x',x_N}f|^2 +
x_N^{2 \xg}|\nabla_{y}f|^2 \right) dx  \right)^{\frac{1}{2}}.\ee
We
 now choose $s=\frac{2}{2-p}$ (so that $sp=2s-2$),  $Q=sp$,  $BQ={\tilde{b}}p$
 and $2A=2 \xb - {\tilde{b}}p$. For this choices we get that
\[
\left(\int_{\R^N_{+}} x_N^{BQ} |f(x)|^Q dx \right)^{\frac{2}{Q}}
\le C \int_{\R^N_{+}} x_N^{2A} \left(|\nabla_{x',x_N} f|^2 +
x_N^{2 \xg}
 |\nabla_{y} f|^2 \right) \ dx\ .
\]
 with  $\xg \in \R$,
$2  \le Q  \le \frac{2N}{N-2}$  if $N\ge 3$ or for $Q\ge 2$ if
$N=2$, and $B=A-1+\frac{Q-2}{2Q}(N+\xg(k-1))$.  The condition $\xb
\neq 0$ is equivalent to $2A+BQ  \neq 0$.

The case where  $f\in C^\infty_0(\R^N)$ and  $2A+BQ > 0$,
 or equivalently $\xb >0$  is practically the same;
we just note that (\ref{4.14}) remains true for $\xb>0$.

\finedim

We are now ready to give the  proof of  Theorem \ref{th9}.
\medskip

\noindent {\em Proof of Theorem \ref{th9}: } We will
 use a (finite) partition of unity for $\Omega$
which we denote by $\varphi_i$, $i=0,\cdots,m$, such that
$1=\sum^m_{i=0} \varphi_i^2$. We  denote by $\Omega_i$ the support
of each  function $\varphi_i$.  We assume  $\Omega_0 \subset
\subset \Omega$ and therefore $c\le d(\lambda,\partial \Omega)\le
c^{-1}$ for  $\lambda \in \Omega_0$.  For  $i \geq 1$, in each
$\xO_i$ we will use  local coordinates  $(y^i,x^{i}_{N})$,  $i\in
\{1,\cdots, m\}$ with  $y^i \in  \Delta_i:=\{y^i: |y^{i}_{j}|\le
\beta \hbox { for } j=1,\cdots, k-1\}$  for some positive constant
$\xb <1$. Each point $\xl \in \bar{\xO}_i  \cap \partial \Omega$
is described by
 $\lambda=(y^i,a_i(y^i))$,
where the functions $a_i$ satisfy a  Lipschitz condition on
$\overline\Delta_i$ with a constant $A>0$ that is
$$|a_i(y^i)-a_i(z^i)|\le A|y^i-z^i|$$ for $y^i, z^i \in
\overline \Delta_i$;  We next define $\hat{B}_i$ by
 $\hat{B}_i:=\{(y^i,x^{i}_{N}): y^i\in
\Delta_i, a_i(y^i)-\beta <x^{i}_{N}<a_i(y^i)+\beta\}$
so that  $\hat{B}_i \cap \Omega = \{(y^i,x^{i}_{N}): y^i\in \Delta_i,
a_i(y^i)-\beta <x^{i}_{N}<a_i(y^i)\}$ and $\Gamma_i=\hat{B}_i \cap
\partial \Omega=\{(y^i,x^{i}_{N}): y^i\in \Delta_i,
x^{i}_{N}=a_i(y^i)\}$. We note that $\xO_i \subset \hat{B}_i \cap \Omega$.
Next we  observe that for any $y\in
 \hat{B}_i \cap \Omega$  we have that  $(1+A)^{-1} (a_i(y^i)-x^{i}_{N}) \le d(\lambda)\le
(a_i(y^i)-x^{i}_{N}), $ (see, e.g.,  Corollary 4.8 in \cite{K}).
By straightening the boundary $\Gamma_i$  we may suppose that
$\Gamma_i \subset \{ x_N^i=0 \}$. From now on we omit the
subscript $i$ for convenience.

As a first step we will prove  that for
 $u\in C^\infty_0(B_1\times H_1^+)$,
 where $B_1:= \{|x'|<1\}$ and
$H_1^+:= \{|y'|<1\} \times  \{0<x_N<1\}$ there holds
\be\la{16.1}
\left(\int_{B_1 \times H_1^+} x_N^\alpha |u|^Q dx'dx_Ndy
\right)^{\frac{2}{Q}} \le C \int_{B_1\times H_1^+} x_N^{\xa}
\left(|\nabla_{x'} u|^2 +x_N^{\sigma}|\nabla_{y} u|^2+
x_N^{\sigma}|\partial_{x_N} u|^2\right) dx'dx_Ndy,
\ee
where $x'\in \R^{N-k}$, and $\lambda =(y, x_N)$ with
 $y\in \R^{k-1}$ and $x_N\in \R$.
 We change
variables by $t=x_N^{\frac{2-\sigma}{2}}$ thus obtaining
\begin{equation} \label{hol} \left(\int_{\{0<t<1\}}
 t^{\frac{2\alpha+\sigma}{2-\sigma}} |u|^Q
dx'dt dy \right)^{\frac{2}{Q}} \le C \int_{\{0<t<1\}}
t^{\frac{2\alpha+\sigma}{2-\sigma}}\left(|\nabla_{x',t} \,  u|^2
+t^{\frac{2\sigma}{2-\sigma}}|\nabla_{y} u|^2\right)dx'dt dy \ .
\end{equation}
Estimate (\ref{hol})  follows from Theorem \ref{th40} part (i) taking
$2A=BQ=\frac{2\alpha+\sigma}{2-\sigma}$ and
$2\gamma=\frac{2\sigma}{2-\sigma}$. We note that $2A+BQ \neq 0$
is equivalent to $2 \xa +\xs \neq 0$. Also, $Q^k_{cr}  \leq \frac{2N}{N-2}$
since $2 \xa + k \xs \geq  0$.

Next, for  $f(x)  = \sum_{i=0}^{m} \varphi_i(\xl)f(x) $, we write
\be\la{16.2} \left(\int_{B_1 \times \xO} d^\alpha |f|^Q dx
\right)^{\frac{2}{Q}} \le C \sum_{i=0}^{m} \left(\int_{B_1 \times
\xO} d^\alpha |f \varphi_i  |^Q dx \right)^{\frac{2}{Q}} \ee Using
(\ref{16.1})  for $i = 1, \ldots, m$ and the standard Sobolev
inequality for $i=0$, in the right hand side of (\ref{16.2}),
after some calculations,
 we end up with
\be\la{16.3}
\left(\int_{B_1 \times \xO} d^\alpha |f|^Q dx
\right)^{\frac{2}{Q}} \le C
   \int_{B_1 \times \xO} d^{\xa}
\left(|\nabla_{x'} f|^2 +  d^{\sigma} |\nabla_{\l} f|^2   \right)
dx + \int_{B_1 \times \xO} d^{\xa+\xs} f^2 dx. \ee To estimate the
last term, we first  use Proposition \ref{lemmata} to obtain
\be\la{16.4} \int_{B_1 \times \xO} d^{\xa+\xs} f^2 dx \leq C
\left( \int_{B_1 \times \xO} d^{\xa +\xs} |\nabla_{\xl} f|^2 dx
+\int_{B_1 \times \xO}  d^{\xa} f^2  dx\right), \ee and then
Poincare inequality in the $x'$--variables,
\[
\int_{B_1 \times \xO} d^{\alpha}f^2 dx \le \int_{B_1 \times \xO}
d^{\xa} |\nabla_{x'} f|^2  dx.
\]
Hence, we end up with
\[
\int_{B_1 \times \xO} d^{\xa+\xs} f^2 dx \leq C
 \int_{B_1 \times \xO} d^{\xa}
\left(|\nabla_{x'} f|^2 +  d^{\sigma} |\nabla_{\l}
f|^2   \right) dx.
\]
Combining this with (\ref{16.3}) we conclude the result.

\finedim

We next prove the Proposition we used in the   proof  of Theorem  \ref{th9}.

\begin{proposition} \label{lemmata}
Let   $\alpha +\sigma<1$.  Then there exists a constant
 $C=C(\xa,\xs, \xO)>0$ such that
\begin{equation}\label{poi2}\int_{\Omega} d^{\alpha+\sigma-2} f^2 d\lambda \le C
\int_{\Omega} d^{\alpha+\sigma} |\nabla_\lambda f|^2 d\lambda \ \
, \ \ \forall \ f\in C^\infty_0(\Omega) ;\end{equation}
The previous inequality fails when $\alpha +\sigma \geq 1$.

Let
$\alpha+\sigma=1$. Then  there exists a constant $C=C(\xa, \xO)>0$
such that
\begin{equation} \label{poi3}\int_{\Omega} \frac{X^2(d)}{d}
f^2 d\lambda \le C \left[\int_{\Omega} d |\nabla_\lambda f|^2
d\lambda +\int_{\Omega} d^\alpha f^2 d\lambda \right]   \ , \ \
\forall \ f\in C^\infty_0(\Omega),
\end{equation}
where $X(d):=(1-\ln(d/D))^{-1}$ and   $D:= \sup_{\lambda \in \xO}d(\lambda)$.

Finally, if $\alpha+\sigma > 1$, there exists a constant $C=C(\xa, \xs, \xO)>0$
such that
\begin{equation}\label{poi1} \int_{\Omega} d^{\alpha+\sigma-2} f^2
d\lambda \le C \left[\int_{\Omega} d^{\alpha+\sigma}
|\nabla_\lambda f|^2 d\lambda +\int_{\Omega} d^\alpha f^2
d\lambda\right] \ , \ \ \forall \ f\in C^\infty_0(\Omega)
.\end{equation}
\end{proposition}

\noindent {\em Proof of Proposition \ref{lemmata}:}

\noindent
{\bf step 1: An auxiliary estimate:} Let
  $\Omega_\delta:=\{\lambda \in \Omega : d(\lambda)\le \delta\}$.
We will establish the following estimate:  Given any $\xe>0$ there
exists a $\xd_0>0$ such that for any $0 <\xd \leq \xd_0$ and any
$u\in C^\infty_0(\Omega)$
\be\la{st1} \int_{\Omega_\delta}
|\nabla_\lambda u|^2 d\lambda \ge \frac{1}{4}\int_{\Omega_\delta}
\frac{u^2}{d^2} d\lambda + \left( \frac{1}{4} - \xe \right)
\int_{\Omega_\delta} \frac{X^2(d) u^2}{d^2} d\lambda
 +\frac{1-X(\xd)}{2 \xd}
\int_{ \partial \Omega_\delta} u^2 dS . \ee To prove this our
starting point is the obvious relation
\[
0 \leq \int_{\Omega_\delta}  \left| \nabla_\lambda u  - \left(
\frac{ \nabla  d}{2 d} - \frac{X  \nabla  d}{2d} \right) u
\right|^2 d \xl.
\]
Expanding the square, integrating by parts and using the fact that
$| d \Delta d |$  can be made arbitrarily small in $\Omega_\delta$,
 for $\xd$ sufficiently small, the result follows.

\noindent
{\bf step 2: Proof of (\ref{poi2})} We change variables by
$u:=d^{\frac{\alpha+\sigma}{2}} f$.   A straightforward calculation
leads to the following identity
\be\la{15.1}
\int_{\Omega_\delta} d^{\alpha+\sigma}|\nabla_{\lambda} f|^2
d\lambda  =\int_{\Omega_\delta} |\nabla_\lambda u|^2 d\lambda
-\frac{\xa+\xs}{2} \left(1-\frac{\xa+\xs}{2}  \right)
\int_{\Omega_\delta} \frac{u^2}{d^2}d\lambda
+\frac{\alpha+\sigma}{2} \int_{\Omega_\delta}
\frac{\Delta d }{d} u^2 d\lambda
-\frac{\alpha+\sigma}{2  \xd }
 \int_{ \partial \Omega_\delta } u^2 dS
\ee From  (\ref{st1}), (\ref{15.1}) and using the fact that $|
\Delta d| <C$ in $\Omega_\delta$, we easily get that there exist
positive constants $c$ such that for $\xd$ sufficiently small
\be\la{15.2} \int_{\Omega_\delta}
d^{\alpha+\sigma}|\nabla_{\lambda} f|^2 d\lambda \geq c
\int_{\Omega_\delta} d^{\alpha+\sigma -2} f^2 d \lambda +
\frac{c}{ \xd^{1- (\alpha+\sigma)} }\int_{ \partial \Omega_\delta
} f^2 dS. \ee On the other hand, away from the boundary we have
that
\[
\int_{\Omega\setminus \Omega_\delta} f^2 \le C
\int_{\Omega\setminus \Omega_\delta} |\nabla_\lambda f|^2 +C
\int_{ \partial \Omega_\delta} f^2 dS,
\]
 from which one can easily deduce
\be\la{15.3}
\int_{\Omega\setminus \Omega_\delta} d^{\alpha+\sigma-2} f^2 \le C_\delta
\int_{\Omega\setminus \Omega_\delta}
d^{\alpha+\sigma}|\nabla_\lambda f|^2 +  \frac{C_\delta}
{\delta^{1-(\sigma+\alpha)}} \int_{ \partial \Omega_\delta } f^2 dS.
\ee
Combining (\ref{15.2}) and (\ref{15.3}) the result follows.

We note that when $\xa+\xs \geq 1$, the constants can be approximated by $C_0^{\infty}(\xO)$
functions in the  norm  given by $\|v\|_{H^{1}(d^{\xa+\xs})}:=\int_{\xO} d^{\xa+\xs}
( |\nabla v|^2 + v^2 ) d \xl$;  see Theorem 2.11 of \cite{FMT}. In particular,
one can  put  a constant function in (\ref{poi2}) to obtain an obvious contradiction.

\noindent {\bf step 3: Proof of (\ref{poi3}) and (\ref{poi1})  }
We first give the proof of (\ref{poi3}).  For $g \in
C_0^{\infty}(\xO_\delta)$ and $u = d^{\frac12} g$ we get  from
(\ref{st1}) and (\ref{15.1}) that for  any $\xe>0$ there exists a
$\xd_0>0$ such that for any $0 <\xd \leq \xd_0$, \be\la{15.4}
\int_{\Omega_\delta} d |\nabla_{\lambda} g|^2 d\lambda    \geq
\left( \frac14 -  2 \xe \right) \int_{\Omega_\delta} \frac{X^2(d)
g^2}{d} d\lambda. \ee To establish the result we argue as follows.
Let $\xi(s)$ be a $C^1$ function such that $0\le \xi\le 1$,
$\xi(s)=0$ if $s \ge 2$, $\xi(s)=1$ if $0\le s\le 1$, and let us
define $\varphi(\lambda)=\xi\left(\frac{d}{\delta}\right)$.

Whence  for  $f\in C^\infty_0(\Omega)$  we have that $f = \varphi
f + (1-\varphi) f$. Using  (\ref{15.4}) for $\varphi f $ and the
fact that on the support of $(1-\varphi) f$ we have
$\frac{X^2(d)}{d}\le \delta^{-1}$ and $d^{\alpha}\ge min
\{\delta^\alpha,D^\alpha\}$, we arrive at \bea \int_{\Omega}
\frac{X^2(d) f^2}{d} d\lambda & \leq & c_1 \int_{\Omega} d
(\varphi^2 |\nabla_\lambda f|^2 + |\nabla_\lambda \varphi|^2 f^2)
d \xl \nonumber + c_2 \int_{\Omega} (1-\varphi)^2 f^2 d \xl
\nonumber \\
& \leq & C \left( \int_{\Omega} d  |\nabla_\lambda f|^2  d \xl  +
\int_{\Omega} d^{\xa} f^2
 d \xl  \right).
\eea
To  prove (\ref{poi1}) we work similarly. We just note that the analogue
of (\ref{15.4}) is
\[
\int_{\Omega_\delta} d^{\xa+\xs} |\nabla_{\lambda} g|^2
d\lambda    \geq  \left( \frac{1-(\xa+\xs)}{2} \right)^2
\int_{\Omega_\delta}d^{\xa+\xs-2} g^2  d\lambda,
~~~~~~~ g \in C_0^{\infty}(\xO_\delta).
\]
We omit further details.

\finedim

\noindent {\bf Remark. }  The  limit  case $k=N$ of Theorem
\ref{th9}, corresponds to the following  isotropic weighted
inequality:
\begin{equation}\label{mava}
\left(\int_{\xO} d^\alpha|f|^Q d\lambda\right)^{\frac{2}{Q}} \le C
\int_{\xO} d^{\xa+\sigma}|\nabla_{\l} f|^2 d\lambda,
\end{equation}
where $\Omega$ is a smooth bounded domain in $\R^N$ and $f\in
C^\infty_0(\Omega)$.  Using similar arguments one can show that
 the above inequality is true if
$\alpha+\sigma<1$,  provided that $\xa >-1$, $2 \xa + N \xs \geq 0$
and  $2 \leq Q \leq Q^N_{cr} = \frac{2(N+\alpha)}{N+\alpha+\sigma-2}$.
On the other hand  inequality  (\ref{mava})  fails
 if $\alpha+\sigma \ge 1$, by the argument of Proposition \ref{lemmata}.

We finally have the following analogue of Corollary \ref{th3}.

\begin{corollary} \label{th10}
For  $N\ge 3$, $1<k<N$, $m>-1$ and $\epsilon\in(0,\frac12)$ we set
\[
C_{1,\epsilon}:=\{(x',\lambda)\in \R^{N-k}\times
\R^k: |x'|<1,~ |\lambda|<1-\epsilon^{1+m}\}\ .
\]
Let $\alpha>-1$ and  $\beta >0$ satisfy
$$-2\alpha(1+m)<\beta m<2(1+m)~ ~~~~{\rm and}~~~~2 \xa (1+m) + \beta k m \geq 0.$$
Then, for any $P$ with
\[
2 \leq P \leq P_{cr}(N,m,\alpha, \beta,k)
 :=  \frac{2 \left( N + \frac{2 \alpha (1+m)+\beta k m}{2 (1+m)-\beta m}\right)  }
{N + \frac{2 \alpha (1+m)+\beta k m}{2 (1+m)-\beta m} -2},
\]
there exists a positive constant $C=C(N,P,m,\alpha,
\beta,k)$ independent of $\epsilon$, such that
 for any function $f\in C^\infty_0(C_{1,\epsilon})$ there holds
$$\left(\int_{C_{1,\epsilon}}
(1-|\lambda|)^\alpha |f(x',\lambda)|^P dx' d\lambda
\right)^{\frac{2}{P}} \le C \int_{C_{1,\epsilon}}(1-|\lambda|)^\alpha
\left(|\nabla_{x'} f|^2
+\frac{(1-|\lambda|)^\beta}{\epsilon^\beta} |\nabla_\lambda
f|^2\right) dx' d\lambda. $$
\end{corollary}

\medskip

\noindent
{\em Proof:} It follows from Theorem \ref{th9}.
 We   have  that  $1-|\lambda|>\epsilon^{1+m}$,
 that is $\epsilon^{-1} >(1-|\lambda|)^{-\frac{1}{1+m}}$ and consequently
$\frac{(1-|\lambda|)^\beta}{\epsilon^\beta}>(1-|\lambda|)^{\frac{\beta
m}{1+m}}$,  $\beta >0$. The result then  follows from Theorem \ref{th9} by
choosing $\sigma:=\frac{\beta m}{1+m}$ there; in particular
 $P_{cr}(N,m,
\alpha, \beta,k) = Q_{cr}(N,\alpha,\frac{\beta m}{1+m},k)$.

\medskip
\noindent {\bf Acknowledgments} The authors are grateful to Prof.
E. Lanconelli for interesting conversations and suggestions which
improved the presentation. The authors would like to thank E.
Cinti for making her thesis available to them after a first draft
of the present paper has been written. LM would like to thank
Prof. X. Cabr\'e, who has drawn her attention to the problem
treated in Sections 2-3, and for several discussions. LM
acknowledges the support of University of Crete and FORTH during
her visit to Greece.  AT acknowledges the support of Universities
of Rome I and Bologna as well as the GNAMPA project "Liouville
theorems in Riemannian and sub-Riemannian settings" during his
visits in Italy.

\end{section}


\end{document}